\documentclass[12pt]{amsart}
\usepackage{amssymb}
\usepackage{amsfonts}

\newtheorem{theorem}{Théorème}
\theoremstyle{plain}

\newtheorem{corollary}{Corollaire}

\newtheorem{definition}{Définition}

\newtheorem{lemma}{Lemme}

\newtheorem{proposition}{Proposition}

\numberwithin{equation}{section}
\input{tcilatex}

\begin{document}
\title[Le Th\'{e}or\`{e}me de Levelt]{Le Th\'{e}or\`{e}me de Levelt}

\begin{abstract}
Let $(E)$ a homogeneous linear differential equation of order $n$ Fuchsien
over $\mathbb{P}^{1}\left( \mathbb{C}\right) $. The idea of Riemann (1857)
was to obtain the properties of solutions of (E) by studying the local
system. Thus, he obtained some properties of Gauss hypergeometric functions
by studying the assocated rank $2$ local system over $\mathbb{P}^{1}\left( 
\mathbb{C}\right) \backslash \left\{ 3\ points\right\} $. For example, he
obtained the Kummer transformations of the hypergeometric functions without
any calculation. The success of the Riemann's methods is due to the fact
that the irreducible rank $2$ local system over $\mathbb{P}^{1}\left( 
\mathbb{C}\right) \backslash \left\{ 3\ points\right\} $ are "rigid". Levelt
theorem, see \cite{B} Theorem 1.2.3 proves this result. In this work, we
propose a partial generalization of this theorem.

\medskip 

\textbf{R\'{E}SUM\'{E}:} Soit $(E)$ une \'{e}quation diff\'{e}rentielle lin%
\'{e}aire homog\`{e}ne Fuchsienne d'ordre $n$ sur $\mathbb{P}^{1}\mathbb{(C})
$. L'id\'{e}e de Riemann (1857) \'{e}tait d'obtenir les propri\'{e}t\'{e}s
des solutions de $(E)$ et ce en \'{e}tudiant le syst\`{e}me local associ\'{e}%
. Ainsi, il obtient certaine propri\'{e}t\'{e}s des fonctions hyperg\'{e}om%
\'{e}triques de Gauss, en \'{e}tudiant le syst\`{e}me locale d'ordre $2$ sur 
$\mathbb{P}^{1}\backslash \left\{ 3\ points\right\} $ associ\'{e}. Par
exemple, il retrouve les transformations de Kummer pour les fonctions hyperg%
\'{e}om\'{e}trique sans faire de calcul. La d\'{e}marche de Riemann a
aboutit car le syst\`{e}me local irr\'{e}ductible associ\'{e} est lin\'{e}%
airement \textquotedblleft rigide\textquotedblleft . Le th\'{e}or\`{e}me de
Levelt, voir \cite{B} Th\'{e}or\`{e}me 1.2.3, prouve ce r\'{e}sultat. Dans
ce travail, on propose une pr\'{e}sentation l\'{e}grement diff\'{e}rente et
une g\'{e}n\'{e}ralisation partielle de ce th\'{e}or\`{e}me.
\end{abstract}

\author{Lotfi Saidane}
\address{Lotfi Saidane, D\'{e}partement de Math\'{e}matiques, Facult\'{e}
des sciences de Tunis, Campus Universitaire, 1060 Tunis, Tunisie.}
\email{lotfi.saidane@fst.rnu.tn}
\thanks{Je remercie le professeur D. Bertrand pour toutes les discussions
concernant ce travail.}
\date{Juin 2009}
\subjclass{12H05}
\keywords{op\'{e}rateur hyperg\'{e}om\'{e}trique, monodromie, systeme lin%
\'{e}airement rigide}
\maketitle

\section{Introduction}

\noindent En langage classique, un module diff\'{e}rentiel sur un corps diff%
\'{e}rentiel $(K,\partial )$ est un $K$-espace vectoriel $M$ de dimension
finie sur $K$ muni d'une application $\nabla _{\partial }:M\rightarrow M$
additive v\'{e}rifiant : 
\begin{equation*}
f\in K,\ m\in M,\ \nabla _{\partial }(fm)=(\partial f)m+f(\nabla _{\partial
}m).
\end{equation*}%
Le fait que la structure de module diff\'{e}rentiel d\'{e}pend de la d\'{e}%
rivation de $K$ conduit \`{a} des difficut\'{e}s dans l'\'{e}tude de
certains probl\`{e}mes de passage du locale au globale, comme par exemple le
cas du probl\`{e}me de Riemann-Hilbert. On a, donc, besoin d'un concept de
module diff\'{e}rentiel beaucoup plus g\'{e}n\'{e}rale, ne d\'{e}pendant pas
de la param\'{e}trisation choisie. Le langage moderne, pour l'\'{e}tude
analytique du probl\`{e}me de Riemann-Hilbert, utilise les concepts de fibr%
\'{e} vectoriel et de syst\`{e}me locale.

\subsection{Connection r\'{e}guli\`{e}re}

\noindent Soit $X$ une surface de Riemann connexe. On d\'{e}signe par $%
\mathcal{O}_{X}$ le sh\'{e}ma des fonctions holomorphes sur $X.$ Un fibr\'{e}
vectoriel $M$ de rang $m$ sur $X$ est un sh\'{e}ma de $\mathcal{O}_{X}-$%
modules sur $X$ localement isomorphe au sh\'{e}ma $\mathcal{O}_{X}^{m}$ de $%
\mathcal{O}_{X}$-modules. On dit que $M$ est un fibr\'{e} trivial s'il est
globalement (sur tout $X$) isomorphe \`{a} $\mathcal{O}_{X}^{m}$. Il est
connu que tout fibr\'{e} vectoriel sur une surface de Riemann connexe non
compact est trivial. On appelle connexion (r\'{e}guli\`{e}re) sur le fibr%
\'{e} vectoriel $M$, tout morphisme de sh\'{e}mas de groupes 
\begin{equation*}
\nabla :M\rightarrow \Omega _{X}\otimes M,
\end{equation*}%
v\'{e}rifiant pour tout ouvert $U$ et tout $f$ $\in \mathcal{O}_{X}(U),$ $%
m\in M(U)$ la relation de Leibniz :%
\begin{equation*}
\nabla (fm)=df\otimes m+f\nabla (m).
\end{equation*}%
Ainsi, si $U$ est un ouvert et $t$ $:$ $U\rightarrow \left\{ c\in \mathbb{C\ 
};\ \left\vert c\right\vert <1\right\} $ est un isomorphisme, alors $\Omega
_{X}(U)$ s'identifie \`{a} $\mathcal{O}_{X}(U)dt$ et $M(U)$ \`{a} $\mathcal{O%
}_{X}^{m}(U).$ L'application%
\begin{equation*}
\nabla (U):M(U)\rightarrow \Omega _{X}(U)dt\otimes M(U)
\end{equation*}%
est une connexion dans le sens classique.

\noindent Soit $X$ un ouvert connexe d'une surface de Riemann $\mathbb{P},$
on suppose que l'infini n'appartient pas \`{a} $X.$ La notion de connexion r%
\'{e}guli\`{e}re $(M,\nabla )$ sur $X$ peut \^{e}tre d\'{e}finie d'une mani%
\`{e}re \'{e}l\'{e}mentaire. On identifie le fibr\'{e} vectoriel $M$ avec $%
\mathcal{O}_{X}^{m},$ le sh\'{e}ma des fonctions holomorphes avec $\mathcal{O%
}_{X}dz.$ Cependant $\nabla $ est determin\'{e} par son action sur $M(X)$,
c'est \`{a} dire par $\nabla _{\frac{d}{dz}}$ sur $M(X).$ Il existe, alors
une matrice $A$ \`{a} coefficients fonctions holomorphes sur $X$ tel que $%
(M,\nabla )$ est compl\`{e}tement d\'{e}termin\'{e} sur $X$ par $\frac{d}{dz}%
+A.$ La connection r\'{e}guli\`{e}re $(M,\nabla )$ sur $X$ est, alors, \'{e}%
quivalente \`{a} $(\mathcal{O}_{X}^{m},\frac{d}{dz}+A).$

\subsection{Fibr\'{e} vectoriel alg\'{e}brique}

\noindent Les fibr\'{e}s vectoriels sur une surface de Riemann ont \'{e}t%
\'{e} classifi\'{e}s par Birkhoff et Grothendieck voir \cite{birkof}, \cite%
{gro}. Pour tout entier $n$, on d\'{e}finit le fibr\'{e} en droite $\mathcal{%
O}_{\mathbb{P}}(n)$ de la mani\`{e}re suivante: soit $U_{0}=\mathbb{P}%
\backslash \left\{ 0\right\} $ et $U_{\infty }=\mathbb{P}\backslash \left\{
\infty \right\} ,$ alors la restriction de $\mathcal{O}_{\mathbb{P}}(n)$ 
\`{a} $U_{0}$ et $U_{\infty }$ est libre et engendr\'{e}e par $e_{0}$ et $%
e_{\infty }$ v\'{e}rifiant sur $U_{0}\cap U_{\infty }$ la relation $%
z^{n}e_{0}=e_{\infty }.$ Un r\'{e}sultat principal montre que tout fibr\'{e}
vectoriel sur la sph\`{e}re de Riemann est isomorphe \`{a} une somme direct
de la forme 
\begin{equation*}
\mathcal{O}_{\mathbb{P}}(n_{1})\oplus ...\oplus \mathcal{O}_{\mathbb{P}%
}(n_{r}),
\end{equation*}%
o\`{u} les $n_{i}$ sont des entiers naturels v\'{e}rifiant $n_{1}\geq
n_{2}\geq ..\geq n_{r}.$ La suite $(n_{1},n_{2},..,n_{r})$ d\'{e}termine le
type du fibr\'{e} $M.$

\noindent D'un point de vue alg\'{e}brique, en associant \`{a} la surface de
Riemann $\mathbb{P}$ la droite projective complexe $\mathbb{P}^{1}(\mathbb{C)%
},$ on peut d\'{e}finir les notions de fibr\'{e}s vectoriels et connections
(alg\'{e}brique), les th\'{e}or\`{e}mes "GAGA" fournissent une \'{e}%
quivalence entre les fibr\'{e}s vectoriels alg\'{e}briques et analytiques.
Soit $X$ un ouvert propre, pour la topologie de Zariski, de la droite
projective$\mathbb{~P}^{1}(\mathbb{C})$. on d\'{e}signe par $\mathcal{O}_{X}$
son alg\`{e}bre affine (le sh\'{e}ma des fonctions r\'{e}guli\`{e}res sur $%
X) $. Soit $M$ un fibr\'{e} vectoriel de rang $m$ sur $\mathbb{P}^{1}.$ La
restriction de $M$ \`{a} $X$ constitute un fibr\'{e} libre, en particulier $%
M(X)$ est un $\mathcal{O}_{X}$-module libre de rang $m.$ La $\mathcal{O}_{X}$%
-alg\`{e}bre des op\'{e}rateurs diff\'{e}rentiels sur $X$, qu'on note $D=$ $%
D_{X}$, est engendr\'{e}e par une d\'{e}rivation $\partial $ de $O_{X}$. On
appelle $D$-module tout $\mathcal{O}_{X}$-module libre de type fini $V$, qui
est muni d'une action $\partial _{V}$ de $\partial $ v\'{e}rifiant les propri%
\'{e}t\'{e}s suivantes:

$\partial_{V}(u+v)=\partial_{V}(u)+\partial_{V}(v)$,

\-$\partial _{V}(fu)=\partial (f)u+f\partial _{V}(u)$, pour tout $u$, $v$
dans $V$ et tout $f$ dans $\mathcal{O}_{X}$.\newline
Dans une base $B$ de $V$ sur $O_{X}$, $\partial _{V}$ est repr\'{e}sent\'{e}
par un op\'{e}rateur $\partial +A_{(B)}$, o\`{u} $A_{(B)}$ est une matrice
carr\'{e} \`{a} co\'{e}fficients dans $O_{X}$. On note $H^{\ast }(V)$ le
groupe de cohomologie du complexe de deRham alg\'{e}bique\ $0\rightarrow
V\rightarrow V\rightarrow 0$ d\'{e}fini par $\partial _{(V)}$. Le groupe $%
H^{0}(V)$ est form\'{e} des sections horizontales (globales), tandis que $%
H^{1}(V)$ s'identifie au conoyau de $\partial _{(V)}$. Pour plus de details,
ainsi que les demonstrations des deux lemmes suivants, voir \cite{Be1} \S .1.

\noindent Soient $V$ et $V^{\prime }$ deux $D$-modules, on note $V^{\ast }$,
le $D$-module $Hom(V,\mathcal{O}_{X})$, de sorte qu'un morphisme de $V$ dans 
$V^{^{\prime }}$est une section horizontale (globale) du $D$-module$\
V^{^{\prime }}\otimes V^{\ast }$. Une extension de $V$ par $V^{^{\prime }}$%
est une suite exacte 
\begin{equation*}
0\rightarrow V^{^{\prime }}\rightarrow E\rightarrow V\rightarrow 0
\end{equation*}%
dans la cat\'{e}gorie des $D$-modules, par abus de langage, on la notera
encore $E$. L'ensemble $Ext_{D}(V,V^{^{\prime }})$ des classes
d'isomorphismes d'extensions de $V$ par $V^{^{\prime }}$ est muni d'une
structure de groupe, qui v\'{e}rifie:

\begin{lemma}
(Coleman): Les groupes $Ext_{D}(V,V^{^{\prime }})$ et $H^{1}(V^{^{\prime
}}\otimes V^{\ast })$ sont canoniquement isomorphe.
\end{lemma}

\noindent Lorsque $V^{^{\prime }}=O_{X}$ (muni de sa connexion canonique).
On d\'{e}duit du lemme pr\'{e}c\'{e}dent que $Ext_{D}(V,O_{X})$ est
canoniquement isomorphe \`{a} $H^{1}(V^{\ast })$ qui a une structure
d'espace vectoriel sur $\mathbb{C}$ (corps de base). Dans le lemme suivant $%
h^{1}(V^{\ast })$ d\'{e}signe la dimension de $H^{1}(V^{\ast })$ et $%
h^{0}(V^{\ast })$ la dimension de $H^{0}(V^{\ast })$ sur $\mathbb{C}$, $S$
est l'ensemble des points a l'infini de $X,$ et $Irr(V)$ la somme des irr%
\'{e}gularit\'{e}s aux diff\'{e}rents points \`{a} l'infini de $S$.

\begin{lemma}
Si $V$ est un $D$-module de rang $n$. Alors: $h^{1}(V^{\ast
})=(cardS-2)n+Irr(V)+h^{0}(V^{\ast })$. En particulier, si $V$ est un module
irr\'{e}ductible non isomorphe \`{a} $O_{X}$ (resp. si $V=O_{X}$ ), il
existe \'{e}xactement $(cardS-2)n+Irr(V)$ (resp. $cardS-1$) extension de $V$
par $O_{X}$ lin\'{e}airement ind\'{e}pendante sur $\mathbb{C}$.
\end{lemma}

\subsection{D-Module}

\noindent Soit $M$ un \'{e}l\'{e}ment unitaire de $D=O_{X}[\partial ]$,
l'action de $\partial $ sur le $O_{X}$ -module $D/DM$, en fait un $D$
-module qu'on note $V(M)$. Les puissances positives ou nulles de $\partial $%
, constitue la base canonique de $V(M)$. Le dual de $V(M)$ est isomorphe 
\`{a} $D/DM^{\ast }$, o\`{u} $M^{\ast }$ d\'{e}signe l'op\'{e}rateur adjoint
de $M$. Les solutions de l'equation diff\'{e}rentielle $M(y)=0$
s'identifient aux vecteurs horizontaux du $D$ -module $V(M^{\ast
})=D/DM^{\ast }$. Par abus d'\'{e}criture on notera de la m\^{e}me mani\`{e}%
re l'equation diff\'{e}rentielle et l'op\'{e}rateur diff\'{e}rentiel associ%
\'{e}.

\noindent Soit $M=L^{\prime }L$ un op\'{e}rateur d\'{e}compos\'{e}, o\`{u} $%
L^{\prime }$ et $L$ sont deux op\'{e}rateurs de $D$ d'ordres strictement
positif. Il d\'{e}finit une extension: 
\begin{equation}
E(L,L^{\prime }):0\longrightarrow V(L^{\prime })\longrightarrow
V(M)\longrightarrow V(L)\longrightarrow 0.  \label{red}
\end{equation}%
Le $O_{X}$-module libre de type fini $V(M)/V(L^{\prime })$ \'{e}tant
isomorphe \`{a} $V(L)$. On choisit la section $s$ du $O_{X}$-module $V(L)$
vers $V(M)$, qui envoie un \'{e}l\'{e}ment $u$ de $V(L)$ sur lui m\^{e}me.
En d\'{e}signant par $B(L)$ et $B(L^{\prime })$ les bases canoniques
respectives de $V(L)$ et $V(L^{\prime })$. L'ensemble $\left\{ B(L^{\prime
})L,B(L)\right\} $ constitue une base de $V(M)$. Si $L=a_{0}+a_{1}\partial
+...+a_{r}\partial ^{r}+\partial ^{r+1}$ et $L^{\prime }=b_{0}+b_{1}\partial
+...+b_{r^{\prime }}\partial ^{r^{\prime }}+\partial ^{r^{\prime }+1}$,
alors l'action de $\partial $ sur $V(L)$, est represent\'{e} par la matrice: 
\begin{equation*}
A_{(L)}=\left( 
\begin{array}{cccc}
0 &  & 0 & -a_{0} \\ 
1 &  & 0 & -a_{1} \\ 
&  &  &  \\ 
0 &  & 1 & -a_{r}%
\end{array}%
\right) ,
\end{equation*}%
l'action de $\partial $ sur $V(L^{\prime })$ est represent\'{e} par la
matrice,

\begin{equation*}
A_{(L^{\prime })}=\left( 
\begin{array}{cccc}
0 &  & 0 & -b_{0} \\ 
1 &  & 0 & -b_{1} \\ 
&  &  &  \\ 
0 &  & 1 & -b_{r^{^{\prime }}}%
\end{array}%
\right) ,
\end{equation*}%
et l'action de $\partial $ sur $V(M)$ est alors repr\'{e}sent\'{e} par la
matrice:%
\begin{equation*}
A_{(M)}=\left( 
\begin{array}{c}
(A_{(L^{\prime })})\left( 
\begin{array}{cc}
0 & 1 \\ 
0 & 0%
\end{array}%
\right) \\ 
\left( 
\begin{array}{cc}
0 & 0 \\ 
0 & 0%
\end{array}%
\right) (A_{(L)})%
\end{array}%
\right) .
\end{equation*}%
Soit $S$ la matrice de $s$ relativement aux bases respectives $B(L)$ et $%
\left\{ B(L^{\prime })L,\text{ }B(L)\right\} $ de $V_{(L)}$ et $V_{(M)}$. La
matrice $S$ est alors \'{e}gale \`{a} $\left( 
\begin{array}{c}
(0) \\ 
(I_{r})%
\end{array}%
\right) $, o\`{u} $I_{r}$ d\'{e}signe la matrice identit\'{e} d'ordre $r$.
On a, par exemple, si $M=(\partial -\alpha )(\partial -\beta ),$ alors $%
A_{(L)}=(\beta ),$ $A_{(L^{\prime })}=(\alpha )$, $A_{(M)}=\left( 
\begin{array}{cc}
\alpha & 1 \\ 
0 & \beta%
\end{array}%
\right) $ et $S=\left( 
\begin{array}{c}
0 \\ 
1%
\end{array}%
\right) .$

\noindent Soit $u$ un \'{e}l\'{e}ment de $V(L)$. On note $(u_{1},...,u_{r})$
ses composantes dans la base $B(L)$. En posant%
\begin{equation*}
\Psi _{(M)}=\partial _{(M)}s-s\partial _{(L)},
\end{equation*}%
le vecteur $\Psi _{(M)}(u)$ est repr\'{e}sent\'{e} par : 
\begin{eqnarray*}
\Psi _{(M)}(u) &=&\partial _{(M)}\left( 
\begin{array}{c}
0 \\ 
. \\ 
0 \\ 
u_{1} \\ 
\\ 
u_{r}%
\end{array}%
\right) -S(\partial _{(L)}(%
\begin{array}{c}
u_{1} \\ 
\\ 
u_{r}%
\end{array}%
)) \\
&=&\left( 
\begin{array}{c}
0 \\ 
. \\ 
0 \\ 
u_{1}^{\prime } \\ 
\\ 
u_{r}^{\prime }%
\end{array}%
\right) +\left( 
\begin{array}{c}
(A_{(L^{\prime })})\left( 
\begin{array}{cc}
0 & 1 \\ 
0 & 0%
\end{array}%
\right) \\ 
\left( 
\begin{array}{cc}
0 & 0 \\ 
0 & 0%
\end{array}%
\right) (A_{(L)})%
\end{array}%
\right) \left( 
\begin{array}{c}
0 \\ 
. \\ 
0 \\ 
u_{1} \\ 
\\ 
u_{r}%
\end{array}%
\right) \\
&&-\left( 
\begin{array}{c}
\left( 0\right) _{r^{\prime }\times r} \\ 
I_{r}%
\end{array}%
\right) \left[ \left( 
\begin{array}{c}
u_{1}^{\prime } \\ 
\\ 
u_{r}^{\prime }%
\end{array}%
\right) +A_{(L^{\prime })}\left( 
\begin{array}{c}
u_{1} \\ 
\\ 
u_{r}%
\end{array}%
\right) \right] \\
&=&\left( 
\begin{array}{c}
u_{r} \\ 
0 \\ 
\\ 
0%
\end{array}%
\right) .
\end{eqnarray*}%
d'o\`{u}%
\begin{equation*}
\Psi _{(M)}(u)=(u_{r},0,..,0).
\end{equation*}%
L'application $\Psi _{(M)}$ est un $O_{X}$-morphisme de $V(L)$ vers $%
V(L^{\prime })$ envoyant le dernier vecteur de $B(L)$ vers le premier
vecteur de $B(L^{\prime }).$ La classe de $\Psi _{(M)}$ dans $%
H^{1}(V_{(L^{\prime })}\otimes V_{(L^{\ast })})$ ne d\'{e}pend pas du choix
de $s$, ni du repr\'{e}sentant de $V_{(M)}$ dans $Ext_{D}(V_{(L)},V_{(L^{%
\prime })}).$

\subsection{Monodromie}

\noindent Soit $S$ un sous-ensemble fini de $\mathbb{C}$, pour $\alpha \in S$%
, soit $A_{\alpha }\in \mathcal{M}_{n}(\mathbb{C)}$, on suppose que $%
\tsum\nolimits_{\alpha \in S}A_{\alpha }=0.$ On pose%
\begin{equation*}
A=\tsum\nolimits_{\alpha \in S}\frac{A_{\alpha }}{z-\alpha }.
\end{equation*}%
On consid\`{e}re le syst\`{e}me diff\'{e}rentiel lin\'{e}aire :%
\begin{equation}
\frac{d}{dz}\left( 
\begin{array}{c}
y_{1} \\ 
\\ 
y_{n}%
\end{array}%
\right) +A\left( 
\begin{array}{c}
y_{1} \\ 
\\ 
y_{n}%
\end{array}%
\right) =0.  \tag{A}
\end{equation}%
Ce syst\`{e}me est singulier r\'{e}gulier sur la droite projective complexe $%
\mathbb{P}^{1}\mathbb{(C})=\mathbb{C}\cup \left\{ \infty \right\} .$ Il
admet des solutions globales multiformes sur $\mathbb{P}^{1}\mathbb{(C}%
)\backslash S.$ L'ensemble de ces solutions constituent un $\mathbb{C}-$%
espace vectoriel $V$ de dimension $n$. Soit $z_{0}\in \mathbb{P}^{1}\mathbb{%
(C})\backslash S,$ l'action du groupe fondamental $\pi _{1}(\mathbb{P}^{1}%
\mathbb{(C})\backslash S,z_{0})$ sur $V$ fournit la repr\'{e}sentaion 
\begin{equation*}
\rho :\pi _{1}(\mathbb{P}^{1}\mathbb{(C})\backslash S,\ z_{0})\rightarrow
GL(V),
\end{equation*}%
appel\'{e}e repr\'{e}sentation de monodromie du syst\`{e}me $(A)$. Du fait
que le groupe $\pi _{1}(\mathbb{P}^{1}\mathbb{(C})\backslash S,\ z_{0})$ est
engendr\'{e} par les classes de lacets $(\gamma _{\alpha })_{\alpha \in S},$
partant de $z_{0}$ puis faisant un tour, dans le sens direct, autour de $%
\alpha ,$ dans un voisinage ne rencontrant pas $S\backslash \left\{ \alpha
\right\} ,$ v\'{e}rifiant la relation (moyennant un ordre convenable de $S)$
:%
\begin{equation*}
\tprod\nolimits_{\alpha \in S}\gamma _{\alpha }=1.
\end{equation*}%
La repr\'{e}sentation $\rho $ est d\'{e}termin\'{e}e par la donn\'{e}e d'un
ensemble de matrices $(M_{\alpha })_{\alpha \in S}$ de $GL_{n}\mathbb{(C})$ v%
\'{e}rifiant 
\begin{equation}
\tprod\nolimits_{\alpha \in S}M_{\alpha }=1  \tag{M}
\end{equation}%
La collection $(M_{\alpha })_{\alpha \in S}$ constitue un syst\`{e}me locale
(complexe) d'ordre $n$ sur $\mathbb{P}^{1}\mathbb{(C})\backslash S$, c'est 
\`{a} dire un sh\'{e}ma d'espaces vectoriels complexe localement isomorphe
au sh\'{e}ma constant $\mathbb{C}^{n}.$ Inversement, soit $(M_{\alpha
})_{\alpha S}$ un ensemble fini de matrices de $GL_{n}\mathbb{(C})$ v\'{e}%
rifiant la relation $\tprod\nolimits_{\alpha \in S}M_{\alpha }=1.$ Le probl%
\`{e}me de Riemann-Hilbert pour les syst\`{e}mes d'\'{e}quations diff\'{e}%
rentielles consiste \`{a} demander si cette collection de matrices peut \^{e}%
tre obtenue de la m\^{e}me mani\`{e}re que pr\'{e}c\'{e}demment, c'est \`{a}
dire la repr\'{e}sentation de monodromie d'un syst\`{e}me diff\'{e}rentiel 
\`{a} p\^{o}le simples en $S$ du type $(A).$

\noindent En 1913, Plemelj et Birkoff \cite{birkof1} r\'{e}alise la repr\'{e}%
sentation $\rho $ comme repr\'{e}sentation de monodromie d'un syst\`{e}me 
\`{a} singularit\'{e}s r\'{e}guli\`{e}res et obtiennent une r\'{e}ponse
positive en supposant que l'un des $M_{\alpha }$ est diagonalisable. En 1928
Lappo-Danileskii \cite{L-D} \'{e}crit explicitement les solutions de $(A)$
sous forme de s\'{e}rie de polylogarithmes, sa r\'{e}ponse est affirmative
dans le cas o\`{u} les matrices $M_{\alpha }$ sont suffisament proche de
l'identit\'{e}. Le point de vue moderne, initi\'{e} par R\"{o}hrl \cite{rohr}%
, consiste \`{a} \'{e}tudier les fibr\'{e}s holomorphes muni d'une connexion
sur une surface de Riemann $X.$ Si $S$ est un sous-ensemble fini non vide de 
$\mathbb{P}^{1}\mathbb{(C})$ et $X=\mathbb{P}^{1}\mathbb{(C})\backslash S,$
alors le fibr\'{e} est trivial, de sorte qu'une section s'identifie \`{a}
une fonction holomorphe de $X$ dans $\mathbb{C}^{n}.$ Il existe, alors, une
application holomorphe $A:X\rightarrow \mathcal{M}_{n}(\mathbb{C})$ telle
que les sections horizontales du fibr\'{e} correspondent aux solutions du
syst\`{e}me diff\'{e}rentiel $Y^{\prime }=AY.$ Les sections horizontales
multiformes constituent un espace vectoriel $V$ de dimension $n,$ sur lequel
le groupe fondamental $\pi _{1}(\mathbb{P}^{1}\mathbb{(C})\backslash S)$ op%
\`{e}re. La repr\'{e}sentation $\rho :\pi _{1}(\mathbb{P}^{1}\mathbb{(C}%
)\backslash S)\rightarrow GL(V)$ est la repr\'{e}sentation de monodromie
associ\'{e}e au fibr\'{e} muni de sa connexion. Le th\'{e}r\`{e}me
d'existence des solutions des syst\`{e}mes d'\'{e}quations diff\'{e}%
rentielles fournit l'\'{e}quivalence entre la cat\'{e}gorie des fibr\'{e}s
holomrphes sur $X$ muni d'une connexion et la cat\'{e}gorie des repr\'{e}%
sentations de $\pi _{1}(X).$ En particulier, toute repr\'{e}sentation de $%
\pi _{1}(\mathbb{P}^{1}\mathbb{(C})\backslash S)$ est la repr\'{e}sentation
de monodromie de $Y^{\prime }=AY,$ o\`{u} $A$ est holomorphe sur $\mathbb{P}%
^{1}\mathbb{(C})\backslash S.$ Deligne \cite{deligne} montre que toute repr%
\'{e}sentation de $\pi _{1}(\mathbb{P}^{1}\mathbb{(C})\backslash S)$ est la
repr\'{e}sentation de monodromie d'un syst\`{e}me $Y^{\prime }=AY$, m\'{e}%
romorphe sur $\mathbb{P}^{1}\mathbb{(C})$ et admettant des singularit\'{e}s r%
\'{e}guli\`{e}res aux points de $S.$ Enfin, Dekkers \cite{dek} montre qu'on
a une r\'{e}ponse positive au probl\`{e}me dans le cas $n=2.$

\noindent La situation est diff\'{e}rente pour les \'{e}quations diff\'{e}%
rentielles lin\'{e}aires scalaires. En effet, soit $n\in \mathbb{N}_{\geq 1}$
et $a_{1},$ .., $a_{n}\in \mathbb{C(}z).$ On consid\`{e}re l'\'{e}quation
diff\'{e}rentielle lin\'{e}aire homog\`{e}ne d'ordre $n$ suivante :%
\begin{equation}
y^{(n)}+a_{1}y^{(n-1)}+....+a_{n-1}y^{\prime }+a_{n}y=0.  \tag{E}
\end{equation}%
On dit que $\alpha \in \mathbb{P}^{1}$ est une singularit\'{e} r\'{e}guli%
\`{e}re de $(E)$ si toute solutions de $(E)$ est \`{a} croissance polyn\^{o}%
miale au voisinage de $\alpha $. Il est connu (voir \cite{Ince}), que $%
\alpha $ est une singularit\'{e} r\'{e}guli\`{e}re de $(E)$ si et seulement
si il existe au moins un $i\in \left\{ 1,..,n\right\} $ tel que $\alpha $
soit un p\^{o}le d'ordre $i$ pour $a_{i}.$ On dit que l'\'{e}quation $(E)$
est Fuchsienne si et seulement si toutes ses singularit\'{e}s sont r\'{e}guli%
\`{e}res. Un calcul simple montre qu'une equation diff\'{e}rentielle lin\'{e}%
aire homog\`{e}ne d'ordre $n$ aynat $s$ singularit\'{e}s r\'{e}guli\`{e}res
dans $\mathbb{P}^{1}$, d\'{e}pend de $\frac{n\left[ n(s-2)+s\right] }{2}$
param\`{e}tres complexes. En revanche, la repr\'{e}sentation de monodromie
de $\pi _{1}(\mathbb{P}^{1}\mathbb{(C})\backslash \left\{ s\ points\right\}
) $ d\'{e}pend de $n^{2}(s-2)+1$ param\`{e}tres. Les deux quantit\'{e}s pr%
\'{e}c\'{e}dentes sont \'{e}gales si et seulement si ($n=1$ et $s$
quelconque) ou ($n=2$ et $s=3).$ Dans ces conditions on peut retrouver l'%
\'{e}quation diff\'{e}rentielle \`{a} partir de la repr\'{e}sentation de
monodromie. Dans le cas g\'{e}n\'{e}ral on ne peut pas obtenir n'importe
quelle repr\'{e}sentation de monodromie comme repr\'{e}sentation de
monodromie d'une \'{e}quation du type de $(E).$ Ainsi, \'{e}tant donn\'{e}
un sous-ensemble fini non vide $S$ de $\mathbb{P}^{1}$ et une collection $%
(M_{\alpha })_{\alpha \in S}$ de matrices de $GL_{n}(\mathbb{C})$ v\'{e}%
rifiant $\tprod\nolimits_{\alpha \in S}$ $M_{\alpha }=1,$ il existe toujours
une \'{e}quation diff\'{e}rentielle lin\'{e}aire r\'{e}guli\`{e}re sur $%
\mathbb{P}^{1}\backslash \Sigma ,$ o\`{u} $\Sigma \supseteq S,$ et singuli%
\`{e}re r\'{e}guli\`{e}re aux points de $S,$ les \'{e}l\'{e}ments \'{e}%
ventuels de $\Sigma \backslash S$ sont des sigularit\'{e}s apparentes pour l'%
\'{e}quation.

\noindent L'id\'{e}e de Riemann (1857) \'{e}tait d'obtenir les propri\'{e}t%
\'{e}s des solutions d'une \'{e}quation diff\'{e}rentielle lin\'{e}aire
fuchsienne d'ordre $n$ par l'\'{e}tude du syst\`{e}me local fournit par le
prolongement analytique au voisinage des diff\'{e}rentes singularit\'{e}s r%
\'{e}guli\`{e}res d'une base de solutions localement holomorphe. Ainsi, il
obtient certaine propri\'{e}t\'{e}s des fonctions hyperg\'{e}om\'{e}triques
de Gauss, en \'{e}tudiant le syst\`{e}me locale d'ordre $2$ sur $\mathbb{P}%
^{1}\backslash \left\{ 3\ points\right\} $ qui leur est associ\'{e}. Il
retrouve les transformations de Kummer pour les fonctions hyperg\'{e}om\'{e}%
trique sans faire de calcul. La d\'{e}marche de Riemann a aboutit car le syst%
\`{e}me local irr\'{e}ductible utilis\'{e} est lin\'{e}airement
\textquotedblleft rigide\textquotedblleft .

\noindent On consid\`{e}re l'\'{e}quation $(E),$ soit $\Sigma $ l'ensemble,
non vide, de ses singularit\'{e}s r\'{e}guli\`{e}re dans $\mathbb{P}^{1}%
\mathbb{(C}).$ On suppose que $(E)$ est r\'{e}guli\`{e}re sur $\mathbb{P}^{1}%
\mathbb{(C})\backslash \Sigma .$ On suppose pour tout point singulier la diff%
\'{e}rence de deux exposants n'est jamais enti\`{e}re. Ce qui se traduit par
le fait que les matrices de monodromie locale poss\`{e}dent des valeurs
propres distinctes. La rigidit\'{e} se traduit par le fait, que si on se
donne une equation diff\'{e}rentielle lin\'{e}aire d'ordre $n$ \`{a}
coefficients dans $\mathbb{C(}z\mathbb{)}$ r\'{e}guli\`{e}re sur $\mathbb{P}%
^{1}\mathbb{(C})\backslash \Sigma $, ayant des singuli\`{e}re r\'{e}guli\`{e}%
re sur $\Sigma $, on suppose que l'\'{e}quation indicielle en tout points de 
$\Sigma $ est la m\^{e}me que celle associ\'{e}e \`{a} $(E),$ alors la
monodromie locale est isomorphe.

\noindent On appelle groupe de monodromie, l'image par la representation de
monodromie du groupe fondamental en un point base $z_{0}\in \mathbb{P}^{1}%
\mathbb{(C})\backslash S$. A conjugaison pr\'{e}s, c'est un sous-groupe bien
defini de $GL_{n}(\mathbb{C}$ $)$. Le groupe de monodromie est un
sous-groupe du groupe de Galois diff\'{e}rentiel. Si les singularit\'{e}s de
l'\'{e}quation sont toute r\'{e}guli\`{e}re, l'adh\'{e}rence de Zariski du
groupe de monodromie remplit le groupe de Galois diff\'{e}rentiel.

\noindent Soit $K$ le corps des fonctions rationnelles sur $X$. On note $%
D_{K}=K[\partial ]$ l'anneau des op\'{e}rateurs diff\'{e}rentiels \`{a}
coefficients dans $K$. Pour les applications que nous avons pr\'{e}vue on
prendra $K=\mathbb{C}(z)$. Tout $D-$ module $V$ definit, par extension des
scalaires un $D_{K}$-module, qu'on notera $V_{K}$, ou encore $V$ (s'il n'y a
pas de risque de confusion).

\section{Op\'{e}rateurs hyperg\'{e}om\'{e}triques}

\noindent Soient $n$ un entier naturel sup\'{e}rieur o\`{u} \'{e}gal \`{a} $%
2 $ et $\alpha _{1}$, ..., $\alpha _{n}$, $\beta _{1}$, ..., $\beta _{n}$
des nombres complexes. On pose $\theta =z\frac{{\Large d}}{{\Large dz}}$, $%
\underline{\alpha }=(\alpha _{1},..,\alpha _{n})$ et $\underline{\beta }%
=(\beta _{1},..,\beta _{n})$. On appelle op\'{e}rateur hyperg\'{e}om\'{e}%
trique (g\'{e}n\'{e}ralis\'{e}) d'ordre $n$ et de param\`{e}tres $\alpha
_{1} $,..,$\alpha _{n}$, $\beta _{1}$,.., $\beta _{n}$ l'op\'{e}rateur :

\begin{equation*}
D(\underline{\alpha },\underline{\beta })=(\theta +\beta _{1}-1)...(\theta
+\beta _{n}-1)-z(\theta +\alpha _{1})...(\theta +\alpha _{n}).
\end{equation*}%
L'\'{e}quation diff\'{e}rentielle $D(\underline{\alpha },\underline{\beta }%
)u=0$ est appel\'{e}e \'{e}quation hyperg\'{e}om\'{e}trique g\'{e}n\'{e}ralis%
\'{e}e d'ordre $n$, elle est r\'{e}guli\`{e}re sur $\mathbb{P}^{1}\mathbb{(C}%
)\backslash \left\{ 0\text{, }1\text{, }\infty \right\} $. Les points $0$, $%
1 $, $\infty $ \'{e}tant des singularit\'{e}s r\'{e}guli\`{e}res pour l'\'{e}%
quation. Les exposants locaux en ces points sont: \newline
$1-\beta _{1}$, ...., $1-\beta _{n}$ en $z=0$; $\alpha _{1}$, ...., $\alpha
_{n}$ en $z=\infty $ et $0,$ $1,$ $..,$ $n-2,$ $-1+\sum_{j=1}^{j=n}(\beta
_{j}-\alpha _{j})$ en $z=1$.

\noindent Si les nombres complexes $\beta _{i}$ sont distincts modulo $%
\mathbb{Z}$. L'equation hyperg\'{e}om\'{e}trique poss\`{e}de $n$ solutions
en z\'{e}ro lin\'{e}airements ind\'{e}pendantes donn\'{e}es par: 
\begin{equation*}
z^{1-\beta _{i}}{}_{n}F_{n-1}(1+\alpha _{1}-\beta _{i},...,1+\alpha
_{n}-\beta _{i};1+\beta _{1}-\beta _{i},...,1+\widetilde{\beta _{i}}-\beta
_{i},..,1+\beta _{n}-\beta _{i}),
\end{equation*}%
pour $i=1$, ..., $n$.

\noindent Le th\'{e}or\`{e}me de Pochhammerr montre que l'\'{e}quation hyperg%
\'{e}om\'{e}trique poss\`{e}de $n-1$ solutions holomorphes lin\'{e}airement
ind\'{e}pendante au voisinage de $z=1.$

\noindent Les r\'{e}sultats suivants, qui existent, en partie dans \cite{BH}%
, nous seront utiles pour la suite de ce paragraphe.

\begin{proposition}
Pour tout $\delta $ dans $\mathbb{C}$ on a:\newline
$(\theta +\delta -1)D(\alpha _{1},...,\alpha _{n};\beta _{1},...,\beta
_{n})= $ $D(\alpha _{1},...,$ $\alpha _{n},$ $\delta ;$ $\beta _{1},...,$ $%
\beta _{n},$ $\delta )$ \newline
et \newline
$D(\alpha _{1},...,$ $\alpha _{n};$ $\beta _{1},...,$ $\beta _{n})(\theta
+\delta )=$ $D(\alpha _{1},...,$ $\alpha _{n},$ $\delta ;$ $\beta _{1},...,$ 
$\beta _{n},$ $\delta +1)$.
\end{proposition}

\begin{proof}[Preuve]
Il suffit de remarquer que pour tout $\delta $ dans $\mathbb{C}$ on a: $%
(\theta +\delta -1)z=z(\theta +\delta )$.
\end{proof}

\begin{corollary}
On a:%
\begin{eqnarray*}
&&D(\alpha _{1},...,\alpha _{n};\beta _{1},...,\beta _{n})(\theta +\alpha
_{j}-1) \\
&=&(\theta +\alpha _{j}-1)D(\alpha _{1},...,\alpha _{j}-1,...,\alpha
_{n};\beta _{1},...,\beta _{n}),
\end{eqnarray*}%
et%
\begin{eqnarray*}
&&D(\alpha _{1},...,\alpha _{n};\beta _{1},...,\beta _{n})(\theta +\beta
_{j}) \\
&=&(\theta +\beta _{j}-1)D(\alpha _{1},...,\alpha _{n};\beta _{1},...,\beta
_{j}+1,...,\beta _{n}).
\end{eqnarray*}
\end{corollary}

\noindent L'op\'{e}rateur $D(\underline{\alpha },\underline{\beta })$ est
dit r\'{e}ductible dans l'anneau $\mathbb{C}(z)[\theta ]$, s'il peut
s'ecrire comme produit de deux op\'{e}rateurs de l'anneau $\mathbb{C}%
(z)[\theta ]$ d'ordre strictement positif. Soient $D_{1}$ et $D_{2}$ deux op%
\'{e}rateurs diff\'{e}rentiels de $\mathbb{C}(z)[\theta ]$ d'ordre
strictement positif, soit $K$ une extension de Picard-Vessiot de $\mathbb{C}%
(z)$ contenant les extensions de Picard-Vessiot correspondantes aux \'{e}%
quations diff\'{e}rentielles $D_{1}(u)=0$ et $D_{2}(u)=0,$ on d\'{e}signe
par $G$ le groupe de Galois diff\'{e}rentiel de $K$. On dira que les op\'{e}%
rateurs $D_{1}$ et $D_{2}$ sont rationnellement \'{e}quivalents si les
espaces des solutions respectifs sont des $G$-modules isomorphes, ou encore,
s'ils existent deux op\'{e}rateurs, $D_{3}$ dans $\mathbb{C}(z)[\theta ]$
sans facteur commun \`{a} droite avec $D_{1}$, et $D_{4}$ dans $\mathbb{C}%
(z)[\theta ]$ tels que: $D_{2}D_{3}=D_{4}D_{1}.$ Autrment dit, l'op\'{e}%
rateur $D_{3}$ r\'{e}alise un isomorphisme entre les espaces des solutions
des \'{e}quations $D_{1}(u)=0$ et $D_{2}(u)=0.$

\noindent Poursuivant \ref{red}, soit $L$ un op\'{e}rateur d'ordre
strictement positif de $\mathbb{C}(z)[\theta ],$ on suppose qu'il s'\'{e}%
crit sous la forme $L=\tprod\nolimits_{i=1}^{r}L_{i}$, o\`{u} les $L_{i}$
sont des op\'{e}rateurs d'ordre strictement positifs irr\'{e}ductibles.
Alors, on obtient une suite de Jordan-H\"{o}lder de module diff\'{e}%
rentielle sur $\mathbb{C}(z)$, telle que chaque facteur est \'{e}quivalent
au module diff\'{e}rentielle associ\'{e}e \`{a} $L_{i}.$ Ainsi, si pour tout 
$i$ on a $L_{i}$ est \'{e}quivalent \`{a} $L_{i}^{\prime }$, alors, $L$ est 
\'{e}quivalent \`{a} $L^{\prime }=\tprod\nolimits_{i=1}^{r}L_{i}^{\prime }$.
L'une des difficult\'{e}s dans la d\'{e}composition d'un op\'{e}rateur, voir %
\ref{Ka}, \ref{boussel}, r\'{e}side dans le fait que si l'un des op\'{e}%
rateurs $L_{i}$ est r\'{e}ductible, alors on ne peut plus continuer le
processus de d\'{e}composition \`{a} \'{e}quivalence pr\'{e}s sans consid%
\'{e}ration de tout l'op\'{e}rateur. On a, par exemple, 
\begin{eqnarray*}
D(0,0,-2;1,1,-1) &=&\theta ^{2}(\theta -2)-z.\theta ^{2}(\theta -2) \\
&=&(1-z).\left[ \theta (\theta -2)\right] .\theta .
\end{eqnarray*}%
L'op\'{e}rateur $\theta (\theta -2)$ est \'{e}quivalent \`{a} $(\theta
-1).(\theta -2)$, en revanche, $D(0,0,-2;1,1,-1)=(1-z).\left[ \theta (\theta
-2)\right] .\theta $ n'est pas \'{e}quivalent \`{a} $\theta (\theta
-1)(\theta -2).$

\noindent Les diff\'{e}rentes d\'{e}compositions d'un op\'{e}rateur hyperg%
\'{e}om\'{e}trique r\'{e}ductibles, fournissent des informations pour le
calcul de son groupe de Galois diff\'{e}rentiel. On a, par exemple, si l'op%
\'{e}rateur hyperg\'{e}om\'{e}trique $L$ s'\'{e}crit%
\begin{equation*}
L=\tprod\nolimits_{i}(\theta +\alpha _{i})L_{L}\tprod\nolimits_{j}(\theta
+\beta _{j}),
\end{equation*}%
o\`{u} $L_{L}$ est Lie-irr\'{e}ductible. Le groupe de Galois diff\'{e}%
rentiel peut-\^{e}tre calcul\'{e} par des extensions succesives de groupes,
voir \ref{boussel} p. 301, elle montre que sous certaine conditions ces
groupes sont ind\'{e}pendants de la d\'{e}composition de l'op\'{e}rateur $L$
et que le groupe $H^{u}=Gal_{diff}(E_{L}/E_{P}E_{L_{L}}E_{Q}),$ o\`{u} $%
P=\tprod\nolimits_{i}(\theta +\alpha _{i}),$ $Q=\tprod\nolimits_{j}(\theta
+\beta _{j}),$ $E_{M}$ est une extension de Picard-V\'{e}ssiot attach\'{e}e 
\`{a} $M$ et $E_{P}E_{L_{L}}E_{Q}$ est un compositum des extensions de
Picard-V\'{e}ssiot, est aussi "gros" que possible. Par exemple, si $L$
(d'ordre $n$) est \'{e}quivalent \`{a} un produit du type $(\theta +\alpha
)(\theta +\alpha ^{\prime })L_{L}$, alors $H^{u}$ est isomorphe \`{a} $%
\mathbb{C}^{2(n-2)}.$

\begin{proposition}
Pour tout entier $s$ dans $\mathbb{Z}$, $D(\alpha _{1}+s,..,\alpha
_{n}+s;\beta _{1}+s,..,\beta _{n}+s)$ est rationnellement \'{e}quivalentsles 
\`{a} $D(\alpha _{1},..,\alpha _{n};\beta _{1},..,\beta _{n})$.
\end{proposition}

\begin{proof}[Preuve]
Soit $s\in \mathbb{Z}$, on a :%
\begin{eqnarray*}
&&D(\alpha _{1},...,\alpha _{n};\beta _{1},...,\beta _{n})z^{s} \\
&=&z^{s}D(\alpha _{1}+s,...,\alpha _{n}+s;\beta _{1}+s,...,\beta _{n}+s)
\end{eqnarray*}%
Du fait que $z^{s}$ est sans facteur en commun \`{a} droite avec $D(\alpha
_{1}+s,...,\alpha _{n}+s;\beta _{1}+s,...,\beta _{n}+s),$ on obtient le r%
\'{e}sultat.
\end{proof}

\begin{proposition}
L'op\'{e}rateur $D(\alpha _{1},...,$ $\alpha _{n};$ $\beta _{1},...,$ $\beta
_{n})$ est r\'{e}ductible si et seulement s'il existe au moins $i$ et $j$
dans \{$1$, $2$, ..., $n$\} tels que $\alpha _{i}-\beta _{j}$ est un entier
relatif.
\end{proposition}

\begin{proof}[Preuve]
Ce r\'{e}sultat est d\'{e}montr\'{e} dans \cite{B}, corollaire 1.2.2, on ne
reprend que la d\'{e}monstration de la condition suffisante. On consid\`{e}%
re les ensembles 
\begin{equation*}
P_{0}=\{(i,j)\ tel\ que\ \beta _{j}-\alpha _{i}\in \mathbb{Z}_{>0}\}
\end{equation*}%
et%
\begin{equation*}
P_{1}=\{(i,j)\ tel\ que\ \beta _{j}-\alpha _{i}\in \mathbb{Z}_{<0}\}.
\end{equation*}%
Si $P_{0}\cup P_{1}=$\{$\emptyset $\}, il existe alors au moins $i$ et $j$
dans \{$1$, ...., $n$\} tel que $\beta _{j}=\alpha _{i}$. Moyennant une
permutation on peut supposer que $\beta _{n}=\alpha _{n}$. De la proposition
1, on d\'{e}duit l'\'{e}galit\'{e}%
\begin{equation*}
D(\underline{\alpha },\underline{\beta })=(\theta +\alpha _{n}-1)D(\alpha
_{1},...,\alpha _{n-1};\beta _{1},...,\beta _{n-1}).
\end{equation*}%
Si $P_{0}\neq $\{$\emptyset $\}, on pose%
\begin{equation*}
n_{0}=Inf_{(i,j)\in P_{0}}\{\beta _{j}-\alpha _{i}\}.
\end{equation*}%
On peut supposer que $\beta _{n}-\alpha _{n}=n_{0}$, on ecrit $\beta
_{n}=\alpha _{n}+1+n_{0}-1$. Si $n_{0}=1$, l'op\'{e}rateur $D$(\underline{$%
\alpha $}, \underline{$\beta $}) est r\'{e}ductible par la proposition 1. Si 
$n_{0}$ $\geq 2$, comme $\beta _{n}-\alpha _{i}\geq n_{0}$ pour tout $i$,
alors $\alpha _{n}+1-\alpha _{i}$ $\notin $ \{ $0$, $-1$, ...... , $2-n_{0}$%
\}. Du lemme 3, on d\'{e}duit que $D(\alpha _{1},...,$ $\alpha _{n};$ $\beta
_{1},...,$ $\alpha _{n}+1+n_{0}-1)$ est rationnellement \'{e}quivalent \`{a} 
$D(\alpha _{1},...,$ $\alpha _{n};$ $\beta _{1},...,$ $\alpha _{n}+1)$ qui
est r\'{e}ductible.\newline
Si $P_{1}\neq $\{$\emptyset $\}. On pose%
\begin{equation*}
n_{1}=sup_{(i,j)\in P_{1}}\{\beta _{j}-\alpha _{i}\}.
\end{equation*}%
Un raisonnement analogue au pr\'{e}c\'{e}dent conduit \`{a} la conclusion.
\end{proof}

\begin{corollary}
Soit $a,$ $b,$ $c\in \mathbb{C}$. Alors l'op\'{e}rateur hyperg\'{e}om\'{e}%
trique%
\begin{equation*}
z(z-1)(\frac{d}{dz})^{2}+\left[ (a+b+1)z-c\right] \frac{d}{dz}+ab
\end{equation*}%
est r\'{e}ductible dans $\mathbb{C}\left( z\right) \left[ \frac{d}{dz}\right]
$ si et seulement si $a$ ou $b$ ou $a-c$ ou $b-c$ est un entier relatif.
\end{corollary}

\begin{proof}
On pose $\theta =z\frac{d}{dz}.$ On a :%
\begin{equation*}
\theta ^{2}=z^{2}\left( \frac{d}{dz}\right) ^{2}+z\frac{d}{dz},
\end{equation*}%
ce qui entra\^{\i}ne%
\begin{eqnarray*}
z\left\{ z(z-1)(\frac{d}{dz})^{2}+\left[ (a+b+1)z-c\right] \frac{d}{dz}%
+ab\right\} &=&\theta \left( \theta +c-1\right) -z\left( \theta +a\right)
\left( \theta +b\right) \\
&=&D(a,b;1,c),
\end{eqnarray*}%
on applique, alors, la proposition pr\'{e}c\'{e}dente pour conclure.
\end{proof}

\begin{proposition}
Tout op\'{e}rateur de $\mathbb{C}(z)[\theta ]$ rationnellement \'{e}%
quivalent \`{a} un op\'{e}rateur r\'{e}ductible est r\'{e}ductibles dans $%
\mathbb{C}(z)[\theta ].$
\end{proposition}

\begin{proof}[Preuve]
Soient $D_{1}$ et $D_{2}$ deux op\'{e}rateurs de $\mathbb{C}(z)[\theta ]$,
de m\^{e}me ordre strictement positif. On d\'{e}signe par $V_{i}$, $i=1$, $2$%
, l'espace des solutions, dans une extension de Picard-Vesiot convenable, de
l'equation diff\'{e}rentielle $D_{i}u=0$.\ Soit $G$ son groupe de Galois diff%
\'{e}retiel, alors l'op\'{e}rateur $D_{1}$ est r\'{e}ductible si et
seulement si il existe un $\mathbb{C}$-espace vectoriel non trivial $G$%
-invariant strictement inclu dans $V_{1}$. Alors l'image, par l'\'{e}%
quivalence rationnelle, de ce sous-espace, est un sous-espace $G$-invariant
de $V_{2}$, ce qui prouve que $D_{2}$ est r\'{e}ductible.
\end{proof}

\bigskip

\noindent Soit $n$ un entier naturel sup\'{e}rieur o\`{u} \'{e}gal \`{a} $2$%
, $\alpha _{1}$, .., $\alpha _{n}$, $\beta _{1}$, .., $\beta _{n}$ des
nombres complexes v\'{e}rifiant l'ensemble $\left\{ \left( i,j\right) \ |\
\alpha _{i}-\beta _{j}\in \mathbb{Z}\right\} $ est non vide. On se propose
de donner quelques propri\'{e}t\'{e}s de l'op\'{e}rateur hyperg\'{e}om\'{e}%
trique $D\left( \alpha _{1},...,\alpha _{n};\beta _{1},...,\beta _{n}\right) 
$. Soit $E=\left\{ \left( i,j\right) \in \left\{ 1,..,n\right\} ^{2}\ |\
\alpha _{i}-\beta _{j}\in \mathbb{Z}\right\} ,$ $E$ est non vide. On pose $%
E=E_{0}\cup E_{+}\cup E_{-},$ o\`{u} $E_{0}=\left\{ \left( i,j\right) \in E\
|\ \alpha _{i}=\beta _{j}\right\} $, $E_{+}=\left\{ \left( i,j\right) \in E\
|\ \alpha _{i}-\beta _{j}\in \mathbb{Z}_{>0}\right\} $ et $E_{-}=\left\{
\left( i,j\right) \in E\ |\ \alpha _{i}-\beta _{j}\in \mathbb{Z}%
_{<0}\right\} $. On pose $t=Inf_{\left( i,j\right) \in E_{-}}\left\{ \beta
_{j}-\alpha _{i}\right\} $ et $s=Sup_{\left( i,j\right) \in E_{+}}\left\{
\beta _{j}-\alpha _{i}\right\} $ On ordonne l'ensemble $\left( \alpha
_{1},..,\alpha _{n};\beta _{1},..,\beta _{n}\right) $ de la mani\`{e}re
suivante :

\noindent On commence par repr\'{e}senter les $\alpha _{i}$ tels que $%
(i,j)\in E_{0},$ ordonn\'{e} par $\mathcal{R}e(\alpha _{i}),$ puis les $%
\alpha _{i}$ tels que $(i,j)\in E_{-}$ et $(i,j)\notin E_{0},$ ordonn\'{e}s
selon la croissance des $\mathcal{\beta }_{j}-\alpha _{i}$, puis les $\alpha
_{i}$ tels que $(i,j)\in E_{+}$ et $(i,j)\notin E_{0}\cup E_{-}$ ordonn\'{e}
par la d\'{e}croissance des $\alpha _{i}-\beta _{j}$ et enfin les $\alpha
_{i}$ tels que $(i,j)\notin E$. Ainsi, on \'{e}crit :%
\begin{equation*}
D\left( \alpha _{1},..,\alpha _{n};\beta _{1},..,\beta _{n}\right) =D(\alpha
_{i\ j})
\end{equation*}%
v\'{e}rifiant :%
\begin{eqnarray*}
\alpha _{1\ 1} &=&\alpha _{1\ 2}=...=\alpha _{1\ r_{1}}=\beta _{1\ 1}=\beta
_{1\ 2}=...=\beta _{1\ s_{1}} \\
&&........................... \\
\alpha _{i_{0}\ 1} &=&\alpha _{i_{0}\ 2}=...=\alpha _{i_{0}\
r_{i_{0}}}=\beta _{j_{0}\ 1}=\beta _{j_{0}\ 2}=...=\beta _{j_{0}\ s_{j_{0}}}
\\
\alpha _{i_{1}\ 1}+t &=&...=\alpha _{i_{1}\ r_{i_{0}}}+t=\beta _{j_{1}\
1}=\beta _{j_{1}\ 2}=...=\beta _{j_{1}\ s_{j_{1}}} \\
&&.........................
\end{eqnarray*}%
Ainsi, voir \ref{boussel} p. 317, s'il existe (exactement) $s$ indices
distincts $i_{1},$ .., $i_{s}$ et (exactement) $s$ indices distincts $j_{1},$
.., $j_{s}$ tels que%
\begin{equation*}
k\in \left\{ 1,..,s\right\} ,\ \ \ \alpha _{i_{k}}-\beta _{j_{k}}\in \mathbb{%
Z}_{\geq 0}
\end{equation*}%
et%
\begin{eqnarray*}
\alpha _{i_{1}}-\beta _{j_{1}} &=&Inf_{(i,j)\in E_{0}\cup E_{+}}(\alpha
_{i}-\beta _{j}) \\
t &\in &\left\{ 2,..,s\right\} ,\ \alpha _{i_{t}}-\beta _{j_{t}}=Inf 
_{\substack{ (i,j)\in E_{0}\cup E_{+}  \\ i\notin \left\{
i_{1},..,i_{t-1}\right\}  \\ j\notin \left\{ j_{1},..,j_{t-1}\right\} }}%
(\alpha _{i}-\beta _{j})
\end{eqnarray*}%
alors il existe $s$ nombres $\alpha _{k}^{\prime }$ tels que $\alpha
_{k}^{\prime }-\beta _{j_{k}}\in \mathbb{Z}$ et l'op\'{e}rateur $D\left(
\alpha _{1},..,\alpha _{n};\beta _{1},..,\beta _{n}\right) $ est \'{e}%
quivalent \`{a}%
\begin{equation*}
\tprod\nolimits_{i\in \left\{ 1,..,s\right\} }(\theta +\alpha _{i}^{\prime
}-1)D(\left\{ \alpha _{i}\right\} _{i\notin \left\{ i_{1},..,i_{s}\right\}
};\left\{ \beta _{i}\right\} _{j\notin \left\{ j_{1},..,j_{s}\right\} }).
\end{equation*}

\noindent En revanche, on a le resultat classique suivant, voir \cite{BH}
corollaire 2.6, est une cons\'{e}quence du lemme pr\'{e}c\'{e}dent.

\begin{proposition}
Si $\beta _{j}-\alpha _{i}$ $\notin $ $\mathbb{Z}$ pour tout indice $i$ et $%
j $ dans \{$1$, ...., $n\}$, alors l'op\'{e}rateur $D$($\alpha _{1}$, ..., $%
\alpha _{n}$; $\beta _{1}$, ..., $\beta _{n}$) est rationnellement \'{e}%
quivalents \`{a} $D$($\alpha _{1}+t_{1}$, ..., $\alpha _{n}+t_{n}$; $\beta
_{1}+s_{1}$, ..., $\beta _{n}+s_{n}$), o\`{u} $t_{1}$, .., $t_{n}$, $s_{1}$,
..., $s_{n}$, sont des entiers relatifs quelconques.
\end{proposition}

\section{Rigidit\'{e} (voir \protect\cite{Katz}, \protect\cite{B})}

\noindent Soit $n\in \mathbb{N}_{\geq 1}$ et $a_{1},$ .., $a_{n}\in \mathbb{%
C(}z).$ On consid\`{e}re l'\'{e}quation diff\'{e}rentielle lin\'{e}aire homog%
\`{e}ne d'ordre $n$ sur $\mathbb{P}^{1}\mathbb{(C})$ suivante :%
\begin{equation}
y^{(n)}+a_{1}y^{(n-1)}+....+a_{n-1}y^{\prime }+a_{n}y=0.  \tag{E}
\end{equation}%
On d\'{e}signe par $S=\left\{ \varpi _{1},..,\varpi _{s}\right\} $
l'ensemble, non vide, de ses singularit\'{e}s (dans $\mathbb{P}^{1}\mathbb{(C%
}))$ qu'on suppose toutes r\'{e}guli\`{e}res. On fixe un point base $%
z_{0}\in \mathbb{P}^{1}\mathbb{(C})\backslash S$ et on d\'{e}signe par $G$
le groupe fondamental $\pi _{1}(\mathbb{P}^{1}\mathbb{(C})\backslash S,\
z_{0}).$ Alors $G$ est un groupe libre engendr\'{e} par les classes
d'homotopies de lacets $\gamma _{i}$ partant de $z_{0},$ faisant un tour
dans le sens direct de $\varpi _{i},$ dans un voisinage ne contenant aucun $%
\varpi _{j},$ $j\neq i,$ puis revenant \`{a} $z_{0},$ tels que 
\begin{equation*}
\tprod\nolimits_{i\in \left\{ 1,..,s\right\} }\gamma _{i}=1.
\end{equation*}%
Du fait que $z_{0}$ est un point r\'{e}gulier pour l'\'{e}quation $(E),$ les
conditions de cauchy sont satisfaites, par cons\'{e}quent, il existe $n$
solutions (locales) de $(E)$ holomorphes au voisinage de $z_{0},$ lin\'{e}%
airement ind\'{e}pendante sur $\mathbb{C}$. Soit $V$ le $\mathbb{C-}$espace
vectoriel constitu\'{e} par ces solutions. La repr\'{e}sentation%
\begin{equation*}
M_{(E)}:\pi _{1}(\mathbb{P}^{1}\mathbb{(C})\backslash S,\ z_{0})\rightarrow
GL(V),
\end{equation*}%
est appel\'{e}e repr\'{e}sentation de monodromie de $(E).$ Pour $i\in
\left\{ 1,..,s\right\} $, on pose $M_{i}=M_{(E)}(\gamma _{i}).$ Alors, $%
M_{i}\in GL_{n}(\mathbb{C})$ et 
\begin{equation*}
\tprod\nolimits_{i\in \left\{ 1,..,s\right\} }M_{i}=I_{n}.
\end{equation*}%
Les matrices $M_{i}$ sont appel\'{e}s matrice de monodromie (locale), elles
constituent un syst\`{e}me locale (complexe) d'ordre $n$ sur $\mathbb{P}^{1}%
\mathbb{(C})\backslash S.$ Le groupe engendr\'{e} par les matrices $M_{i}$
est appel\'{e} "groupe de monodromie" de $(E),$ li\'{e} \`{a} la base de de
solutions locales en $z_{0}.$

\noindent \textbf{Question : }Si on change la base de solution locale, alors
le noiveau groupe de monodromie est-il isomorphe \`{a} l'ancien? En d'autres
termes : Le syst\`{e}me locale d\'{e}finie par les $M_{i}$ est-il "lin\'{e}%
airement rigide"?

\noindent Soient $r\in \mathbb{N}_{\geq 2}$ et $g_{1}$, $g_{2}$, $....$, $%
g_{r}$ des \'{e}l\'{e}ments de $GL_{n}(\mathbb{C})$ v\'{e}rifiant%
\begin{equation*}
g_{1}.g_{2}....g_{r}=Id_{n}
\end{equation*}%
On dit que le $r$-uplet $\left\{ g_{1}\text{, }g_{2}\text{, }...\text{, }%
g_{r}\right\} $ est lin\'{e}airement rigide, si pour tout conjugu\'{e}s $%
\widetilde{g_{1}}$, $\widetilde{g_{2}}$, $...$, $\widetilde{g_{r}}$ de $%
g_{1} $, $g_{2}$, $....$, $g_{r}$ dans $GL_{n}(\mathbb{C}$ $)$ v\'{e}%
rifiant: 
\begin{equation*}
\breve{g}_{1}\breve{g}_{2}...\breve{g}_{r}=Id_{n},
\end{equation*}%
il existe $u$ dans $GL_{n}(\mathbb{C}$ $)$ tel que $\breve{g}%
_{i}=ug_{i}u^{-1}$ pour $i=1,2,...r$. A titre d'exemple, le couple $(g$, $%
g^{-1}),$ $g\in GL_{n}(\mathbb{C}),$ est lin\'{e}airement rigide.

\noindent Le groupe $<g_{1},..,g_{r}>$ est irr\'{e}ductible si et seulement
il agit irr\'{e}ductiblement sur $\mathbb{C}^{n}.$ Katz \cite{Katz} th\'{e}or%
\`{e}me 1.1.2, a crast\'{e}ris\'{e} les formes normales de Jordan des syst%
\`{e}mes locaux irr\'{e}ductibles lin\'{e}airement rigides. Le th\'{e}or\`{e}%
me de Levelt, \cite{B} th\'{e}or\`{e}me 1.2.3, montre que le syst\`{e}me
local associ\'{e} \`{a} une \'{e}quation hyperg\'{e}om\'{e}trique irr\'{e}%
ductible est lin\'{e}airement rigide. Le but de ce paragraphe est de g\'{e}n%
\'{e}raliser en partie le th\'{e}or\`{e}me de Levelt.

\begin{definition}
On dit que $h\in GL(n,\mathbb{C)}$ est une pseudo-reflection si le rang de $%
(h-Id_{n})$ est \'{e}gal \`{a} $1.$
\end{definition}

\begin{lemma}
Soit $n,p\in \mathbb{N}_{\geq 2}$ et $A_{1},$ $A_{2},$ ..., $A_{p}\in $ $%
GL_{n}(\mathbb{C)}$ tels que pour tout $i,$ $j\in \left\{ 1,2,...,p\right\} $%
, $i<j$, on ait $A_{i}A_{j}^{-1}$ est une pseudo-reflection. Alors $A_{1},$ $%
A_{2},$ ..., $A_{p}$ poss\`{e}dent $(n-1)$ lignes ou colonnes communes.\label%
{ligne}
\end{lemma}

\begin{proof}
On fera une d\'{e}monstration pour $p=3,$ le cas g\'{e}n\'{e}ral \'{e}tant
similaire. On suppose que $n\geq 3.$\newline
On pose $W_{1}=\ker (A_{1}-A_{2})$ et $W_{2}=\ker (A_{2}-A_{3}).$ Du fait
que $A_{1}A_{2}^{-1}$ et $A_{2}A_{3}^{-1}$ sont des pseudo-reflections, on d%
\'{e}duit que $W_{1}$ et $W_{2}$ sont des sous-espaces vectoriels de $%
\mathbb{C}^{n}$ de dimension $n-1$. Si $W_{1}=W_{2},$ on choisit une base de 
$W_{1}$ qu'on compl\`{e}te, par un vecteur, en une base de l'espace total.
Relativement \`{a} cette base les matrices $A_{1}$, $A_{2}$ et $A_{3}$ ont
les m\^{e}me $(n-1)$ premi\`{e}re colonnes.\newline
On suppose que $W_{1}\neq W_{2}.$ Du fait que $n\geq 3,$ l'espace vectoriel $%
W_{1}\cap W_{2}$ est de dimension $n-2.$ Soit $\left\{
e_{1},...,e_{n-2}\right\} $ une base de l'espace $W_{1}\cap W_{2}$.\newline
* Si $A_{1}A_{2}^{-1}$ est une reflection, il existe, alors, $e_{n}$ tel que 
$(A_{1}-A_{2})(e_{n})=\nu e_{n},$ avec $\nu \neq 0.$ On choisit $e_{n-1}$
dans $W_{1}$ de sorte que $\left\{ e_{1},...,e_{n-1}\right\} $ soit une base
de $W_{1}.$ Le syst\`{e}me $\underline{e}=\left\{ e_{1},...,e_{n}\right\} $
constitue, alors, une base de $\mathbb{C}^{n}.$ Pour $i\in \left\{
1,2,3\right\} ,$ $j\in \left\{ 1,2,..,n\right\} ,$ On d\'{e}signe par $%
A_{i,j}$, la $j-\grave{e}me$ colonne, par rapport \`{a} la base $\underline{e%
},$ de la matrice $A_{i}$. Du fait que $%
rang(A_{1}-A_{3})=rang(A_{2}-A_{3})=1 $ et $W_{1}\neq W_{2},$ ils existent $%
\lambda $ et $\beta $ non nuls dans $\mathbb{C}$ tels que 
\begin{equation*}
A_{1,n-1}-A_{3,n-1}=\lambda (A_{1,n}-A_{3,n})
\end{equation*}%
et%
\begin{equation*}
A_{2,n-1}-A_{3,n-1}=\beta (A_{2,n}-A_{3,n}).
\end{equation*}%
On a, par hypoth\`{e}se $A_{1,n-1}=A_{2,n-1}$ et $A_{1}e_{n}=A_{2}e_{n}+\nu
e_{n}.$ En rempla\c{c}ant ces relations dans les \'{e}galt\'{e}s pr\'{e}c%
\'{e}dentes, on obtient, par abus d' \'{e}criture,%
\begin{equation*}
\lambda (A_{2,n}+\nu e_{n}-A_{3,n})=\beta (A_{2,n}-A_{3,n}).
\end{equation*}%
Par cons\'{e}quent, il existe au moin $\alpha \in \mathbb{C}$, $\alpha \neq
0 $ tel que 
\begin{equation*}
A_{2,n}-A_{3,n}=\alpha e_{n},
\end{equation*}%
ce qui entra\^{\i}ne%
\begin{equation*}
A_{1,n-1}-A_{3,n-1}=\lambda \alpha e_{n},
\end{equation*}%
et prouve que les $\left( n-1\right) $ premi\`{e}res ligne des matrices $%
A_{1},$ $A_{2}$ et $A_{3}$ sont identiques.\newline
* Si $A_{1}A_{2}^{-1}$ est idempotente (n'admet que $1$ comme valeur
propre), alors, l'image de $A_{1}-A_{2}$ est contenue dans son noyau $W_{1}$%
. Soit $w$ un g\'{e}n\'{e}rateur de $\func{Im}\left( A_{1}-A_{2}\right) ,$
il existe, alors, un vecteur $e_{n}$ dans $\mathbb{C}^{n}$ tel que $\left(
A_{1}-A_{2}\right) e_{n}=w.$ Si $w\in W_{1}\cap W_{2},$ on pose $e_{1}=w$ de
sorte que $\left\{ e_{1},...,e_{n-2}\right\} $ soit une base de l'espace $%
W_{1}\cap W_{2}.$ On choisit $e_{n-1}$ dans $W_{1}$ tel que $\left\{
e_{1},...,e_{n-1}\right\} $ soit une base de $W_{1}.$ Le syst\`{e}me $%
\underline{e}=\left\{ e_{1},...,e_{n}\right\} $ constitue, alors, une base
de $\mathbb{C}^{n}.$ Si $w\notin W_{1}\cap W_{2},$ on pose $e_{n-1}=w.$
Ainsi, il existe un unique indice $m\in \left\{ 1,n-1\right\} ,$ tel que $%
w=e_{m}.$ Dans ces conditions, les $(n-1)$ premi\`{e}re colonnes de $%
A_{1}-A_{2}$ sont nulles et la derni\`{e}re colonne est \'{e}gale \`{a} $%
e_{m}.$ Les $(n-2)$ premi\`{e}res colonnes de $A_{1}-A_{3}$ et de $%
A_{2}-A_{3}$ sont nulles. En \'{e}crivant que $A_{1}-A_{3}$ et $A_{2}-A_{3}$
sont de rang $1,$ on obtient que toutes les composantes de leurs deux derni%
\`{e}res colonnes, sauf la $m-i\grave{e}me$ ligne sont nulles. Par cons\'{e}%
quent, \`{a} l'exeption de la $m-i\grave{e}me$ ligne, toutes les lignes de $%
A_{1}$, $A_{2}$ et $A_{3}$ sont identiques.\newline
Si $n=2$ et $W_{1}\neq W_{2}.$ On refait le m\^{e}me raisonnement en
remplacant $W_{1}\cap W_{2}$ par $\left\{ 0\right\} .$
\end{proof}

\noindent Les deux resultats suivants modifient l'\'{e}nonc\'{e} et g\'{e}n%
\'{e}ralisent, en partie, le th\'{e}or\`{e}me 1.2.1 de \cite{B}.

\begin{theorem}
Soit $n,p\in \mathbb{N}_{\geq 2}$ et $A_{1},$ $A_{2},$ ..., $A_{p}\in $ $%
M_{n}(\mathbb{C)}$ ayant $\left( n-1\right) $ colonnes (ou lignes) communes
et une valeur propre en commun. Alors ces matrices stabilisent au moins une
droite ou un hyperplan de $\mathbb{C}^{n}.\label{stabilise}$
\end{theorem}

\begin{proof}
On remarque qu'un syst\`{e}me de matrices stabilisent un m\^{e}me hyperplan
si et seulement si leurs transpos\'{e}es, comme endomorphisme de l'espace
dual, stabilisent une m\^{e}me droite, autrement dit ils poss\`{e}dent un
vecteur propre en commun. Sans perte de g\'{e}n\'{e}ralit\'{e}s, on pourra
supposer que les $(n-1)$ premi\`{e}res lignes des matrices $A_{i}$ sont
identiques. Si $\lambda $ d\'{e}signe la valeur propre commune, alors les
matrices $A_{1}-\lambda I_{n}$, $A_{2}-\lambda I_{n},$ ..., $A_{p}-\lambda
I_{n}$ ont les m\^{e}me $(n-1)$ premi\`{e}res lignes et sont de rang inf\'{e}%
rieur ou \'{e}gal \`{a} $n-1$ :\newline
- Si les $(n-1)$ premi\`{e}res lignes de ces derni\`{e}res matrices sont lin%
\'{e}airement ind\'{e}pendantes, alors la derni\`{e}re ligne de chacune de
ces matrices sera combinaison lin\'{e}aire des pr\'{e}c\'{e}dentes. Soit $v$
un vecteur non nul orthogonal \`{a} ces $(n-1)$ premi\`{e}re lignes. Alors $%
v $ est orthogonal \`{a} la derni\`{e}re ligne de chacune de ces matrices.
Ainsi, $v$ est un vecteur propre commun, de valeur propre $\lambda ,$ des
matrices $A_{1},$ .., $A_{p}.$\newline
- Si les $(n-1)$ premi\`{e}res lignes de ces matrices sont lin\'{e}airement d%
\'{e}pendantes. Alors leurs transpos\'{e}es ont $(n-1)$ colonnes communes lin%
\'{e}airement d\'{e}pendantes. Soit $c_{1},$ $c_{2},$ ..., $c_{n-1}$ les
coefficients d'une relation de d\'{e}pendance lin\'{e}aire, non triviale,
entre ces colonnes. Alors, le vecteur $v=(c_{1},$ $c_{2},$ ..., $c_{n-1},$ $%
0)$ est un vecteur propre commun de toutes les matrices $A_{1}^{T},$ .., $%
A_{p}^{T}$ de valeur propre $\lambda .$ Par cons\'{e}quent, les matrices $%
A_{1},$ .., $A_{p}$ stabilisent simultan\'{e}ment un hyperplan.\newline
Si les matrices $A_{1},$ .., $A_{p}$ poss\`{e}dent $(n-1)$ colonnes
communes, leus transpos\'{e}es poss\`{e}dent, alors, $(n-1)$ lignes communes
et le raisonnement pr\'{e}c\'{e}dent conduit \`{a} la conclusion.
\end{proof}

\begin{theorem}
Soit $n,p\in \mathbb{N}_{\geq 2}$, $A_{1},$ $A_{2},$ ..., $A_{p}\in $ $M_{n}(%
\mathbb{C)}$ ayant $\left( n-1\right) $ lignes (ou colonnes) communes et
stabilisent un m\^{e}me sous-espace non trivial de $\mathbb{C}^{n}.$ Alors, $%
\cap _{i=1}^{p}specA_{i}\neq \emptyset .\label{spec}$
\end{theorem}

\begin{proof}
On peut supposer que, relativement \`{a} une certaine base $B=\left\{
e_{1},..,e_{n}\right\} $ de $\mathbb{C}^{n},$ les matrices $A_{1},$ ..., $%
A_{p}$ ont les m\^{e}me $(n-1)$ premi\`{e}res colonnes. On d\'{e}signe par $%
E $ le sous-espace de $\mathbb{C}^{n}$ engendr\'{e} par $\left\{
e_{1},..,e_{n-1}\right\} .$ Soit $W$ un sous-espace vectoriel non trivial de 
$\mathbb{C}^{n}$ stable sous l'action des $A_{i}.$ On suppose que $W\subset
E $ et que $\dim _{\mathbb{C}}W=r\in \left\{ 1,..,n-1\right\} .$ Soit $%
\left\{ w_{1},..,w_{r}\right\} $ une base de $W,$ qu'on compl\`{e}te de
sorte que $B^{\prime }=\left\{ w_{1},..,w_{n}\right\} $ soit une base de $%
\mathbb{C}^{n}.$ Du fait que $A_{i}e_{j}=A_{k}e_{j},$ pour $i,k\in \left\{
1,..,p\right\} $ et $j\in \left\{ 1,..,n-1\right\} ,$ on d\'{e}duit que 
\begin{equation*}
A_{i}w_{j}=A_{k}w_{j},~pour~i,k\in \left\{ 1,..,p\right\} \ et\ j\in \left\{
1,..,r\right\} ,
\end{equation*}%
ce qui prouve, compte tenue de la stabilit\'{e} de $W$, que relativement 
\`{a} la base $B^{\prime }$ les matrices $A_{i}$ sont de la forme 
\begin{equation*}
A_{i}=\left( 
\begin{array}{cc}
A_{(r,r)} & \ast _{(r,n-r)}\ \ \ \  \\ 
o & \ast _{(n-r,n-r)}%
\end{array}%
\right) ,
\end{equation*}%
o\`{u} $A_{(r,r)}$ est une matrice d'ordre $r$ commune \`{a} tous les $%
A_{i}, $ $\ast _{(r,n-r)}$ (resp. $\ast _{(n-r,n-r)})$ est un \'{e}l\'{e}%
ment de $M_{(r,n-r)}(\mathbb{C})$ (resp. $M_{(n-r,n-r)}(\mathbb{C})).$ Par
cons\'{e}quent, le polyn\^{o}me $\det (A_{(r,r)}-\lambda I_{r})$ divise tous
les polyn\^{o}mes caract\'{e}ristiques de tous les $A_{i}$. Les racines
complexes de $\det (A_{(r,r)}-\lambda I_{r})$ sont dans $\cap
_{i=1}^{p}specA_{i}.$\newline
On suppose que $W\nsubseteq E,$ il existe, alors, une base $\left\{
f_{n-p+1},..,f_{n}\right\} $ de $W$ et un syst\`{e}me libre $\left\{
g_{1},..,g_{n-p}\right\} $ de $E$ tel que $\left\{
g_{1},..,g_{n-p},f_{n-p+1},..,f_{n}\right\} $ soit une base de $\mathbb{C}%
^{n}.$ Relativement \`{a} cette base les matrices $A_{i}$ sont de la forme :%
\begin{equation*}
A_{i}=\left( 
\begin{array}{c}
A_{(n,n-p)}%
\end{array}%
\begin{array}{c}
0\ \ \ \ \ \  \\ 
\ \ast _{(p,p)}%
\end{array}%
\right) ,
\end{equation*}%
o\`{u} $A_{(n,n-p)}$ est une matrice ayant $n$ lignes et $(n-p)$ colonnes
commune \`{a} tous les $A_{i}$ et $\ast _{(p,p)}$ est une matrice d'ordre $%
p. $ Ainsi, on d\'{e}duit que les $A_{i}$ poss\`{e}dent au moin une valeur
propre commune.
\end{proof}

\begin{corollary}[Beukers Th\'{e}or\`{e}me 1.2.1]
Soit $H$ un sous-groupe de $GL_{n}(\mathbb{C)}$ engendr\'{e} par deux
matrices $A$ et $B$ v\'{e}rifiant $AB^{-1}$ est une pseudo-reflection.
Alors, $H$ est lin\'{e}airement irr\'{e}ductible si et seulement si, les
spectres de $A$ et de $B$ sont disjoints.
\end{corollary}

\begin{proof}
Le lemme \ref{ligne} pour $p=2,$ montre que $A$ et $B$ poss\`{e}dent $(n-1)$
lignes ou colonnes communes. Les th\'{e}or\`{e}mes \ref{stabilise} et \ref%
{spec} permettent de conclure.
\end{proof}

\noindent Le r\'{e}sultat suivant modifie l\'{e}g\`{e}rement et g\'{e}n\'{e}%
ralise, en partie, le th\'{e}or\`{e}me de Levelt au cas $p\geq 2$ (voir \cite%
{B} th\'{e}or\`{e}me 1.2.3).

\begin{theorem}
Soit $n,p\in \mathbb{N}_{\geq 2}$ et $\alpha _{i}=\left\{ \alpha
_{i,1},...,\alpha _{i,n}\right\} \subset \mathbb{C}^{\ast },$ $1\leq i\leq
p, $ v\'{e}rifiant $\cap _{i=1}^{p}\alpha _{i}=\emptyset .$ Alors il existe $%
A_{1},$ $A_{2},$ ..., $A_{p}$ dans $GL_{n}(\mathbb{C)}$ ayant $(n-1)$
colonnes en commun, unique \`{a} conjuguaison pr\'{e}s par un m\^{e}me
isomorphisme, v\'{e}rifiant, pour tout $i,$ $specA_{i}=\alpha _{i}.\label%
{lev}$
\end{theorem}

\begin{proof}
Existence :\newline
Pour tout $(i,j)\in \left\{ 1,..,p\right\} \times \left\{ 1,..,n\right\} ,$
on d\'{e}finit $A_{i,j}$ par 
\begin{equation*}
\tprod\nolimits_{j=1}^{n}(X-\alpha
_{i,j})=X^{n}+\tsum\nolimits_{k=0}^{n-1}A_{i,n-k}X^{k},\ 
\end{equation*}%
du fait que les $\alpha _{i,j}$ sont non nuls, les matrices d\'{e}finient
par 
\begin{equation*}
A_{i}=\left( 
\begin{array}{ccc}
0\ 0 & 0 & -A_{i,n} \\ 
1\ 0 &  & . \\ 
0\ 0 & 1 & -A_{i,1}%
\end{array}%
\right) ,
\end{equation*}%
sont dans $GL(n,\mathbb{C)}$ de polyn\^{o}me caract\'{e}ristique 
\begin{eqnarray*}
\det (XI_{n}-A_{i}) &=&X^{n}+\tsum\nolimits_{k=0}^{n-1}A_{i,n-k}X^{k} \\
&=&\tprod\nolimits_{j=1}^{n}(X-\alpha _{i,j}).
\end{eqnarray*}%
ce qui prouve l'existence.\newline
Unicit\'{e} :\newline
Soit $A_{1},$ $A_{2},$ ..., $A_{p}\in GL(n,\mathbb{C)}$ ayant $(n-1)$
colonnes en commu. On pourra supposer que les $(n-1)$ prem\`{e}re colonnes
des matrices $A_{i}$ sont identiques. Soit $\left\{ e_{1},..,e_{n}\right\} $
une base de $\mathbb{C}^{n}$ par rapport \`{a} laquelle on a, pour $i,j\in
\left\{ 1,..,p\right\} ,$ $k\in \left\{ 1,..,n-1\right\} $ :%
\begin{equation*}
A_{i}e_{k}=A_{j}e_{k}.
\end{equation*}%
On d\'{e}signe par $W$ le sous-espace vectoriel de $\mathbb{C}^{n}$ engendr%
\'{e} par $\left\{ e_{1},..,e_{n-1}\right\} .$ On a :%
\begin{equation*}
\dim _{\mathbb{C}}(W\cap A_{1}W\cap ...\cap A_{1}^{n-2}W)\geq 1.
\end{equation*}%
On suppose que la dimension de $W\cap A_{1}W\cap ...\cap A_{1}^{n-2}W$ est
sup\'{e}rieur ou \'{e}gal \`{a} $2.$ Par cons\'{e}quent, on a 
\begin{equation*}
\dim _{\mathbb{C}}(W\cap A_{1}W\cap ...\cap A_{1}^{n-1}W)\geq 1.
\end{equation*}%
D'o\`{u}, le sous-espace $W\cap A_{1}W\cap ...\cap A_{1}^{n-1}W$ de $\mathbb{%
C}^{n}$ est non trivial et stable sous l'action de tous les $A_{i}.$ Le th%
\'{e}or\`{e}me \ref{spec} montre que $\cap _{i=1}^{p}specA_{i}\neq \emptyset 
$ ce qui est absurde. Ainsi,%
\begin{equation*}
\dim _{\mathbb{C}}(W\cap A_{1}W\cap ...\cap A_{1}^{n-2}W)=1,
\end{equation*}%
il existe, alors, un vecteur $v$ dans $W$ tel que le syst\`{e}me $\left\{
v,A_{1}v,..,A_{1}^{n-2}v\right\} $ constitue une base de $W,$ qu'on compl%
\`{e}te par un vecteur en une base de l'espace total $\mathbb{C}^{n}.$
Relativement \`{a} cette derni\`{e}re les matrices $A_{i}$ sont de la forme%
\begin{equation*}
A_{i}=\left( 
\begin{array}{ccc}
0\ 0 & 0 & -A_{i,n} \\ 
1\ 0 &  & . \\ 
0\ 0 & 1 & -A_{i,1}%
\end{array}%
\right) ,
\end{equation*}%
o\`{u} les $A_{i,j}$ sont determin\'{e}s par le spectre $\left\{ \alpha
_{i,1},...,\alpha _{i,n}\right\} $ de $A_{i}$ de la mani\`{e}re suivante :%
\begin{equation*}
\tprod\nolimits_{j=1}^{n}(X-\alpha
_{i,j})=X^{n}+\tsum\nolimits_{k=0}^{n-1}A_{i,n-k}X^{k},
\end{equation*}%
ce qui termine la preuve.
\end{proof}

\end{document}